\documentclass[11pt]{article}

\usepackage{amsmath,amssymb,amsthm}

\newtheorem{theorem}{Theorem}
\newtheorem{definition}[theorem]{Definition}
\newtheorem{proposition}[theorem]{Proposition}
\newtheorem{remark}[theorem]{Remark}

\def\Tr{\mathrm{Tr}}
\def\LL{{\cal{L}}}
\def\ot{\otimes}
\def\ov{\overline}
\def\un{\underline}
\def\la{{\lambda}}

\def\H{\mathcal{H}}
\def\S{\mathcal{S}}
\def\C{\mathbb{C}}

\def\be{\begin{equation}}
\def\ee{\end{equation}}

\begin{document}

\makeatletter
\renewcommand{\theequation}{{\thesection}.{\arabic{equation}}}
\@addtoreset{equation}{section} \makeatother

\title{Quantum Schur-Weyl duality and $q$-Frobenius formula related to Reflection Equation algebras}
\author{\rule{0pt}{7mm} Dimitry Gurevich\thanks{gurevich@ihes.fr}\\
{\small\it Institute for Information Transmission Problems}\\
{\small\it Bolshoy Karetny per. 19,  Moscow 127051, Russian Federation}\\
\rule{0pt}{7mm} Pavel Saponov\thanks{Pavel.Saponov@ihep.ru}\\
{\small\it
National Research University Higher School of Economics,}\\
{\small\it 20 Myasnitskaya Ulitsa, Moscow 101000, Russian Federation}\\
{\small \it and}\\
{\small \it
Institute for High Energy Physics, NRC "Kurchatov Institute"}\\
{\small \it Protvino 142281, Russian Federation}}

\maketitle

\begin{abstract}
We establish a $q$-version of the Schur-Weyl duality, in which the role of the symmetric group is played by the Hecke algebra and the role of the enveloping algebra $U(gl(N))$ is
played by the Reflection Equation algebra, associated with any skew-invertible Hecke symmetry. Also, in each Reflection Equation algebra we define
analogues  of the Schur
polynomials and power sums in two forms: as polynomials in generators of a given Reflection Equation algebra and in terms of the so-called eigenvalues of the generating matrix $L$, defined by means of the Cayley-Hamilton identity. It is shown that on any Reflection Equation algebra there exists a formula, which brings into correlation the Schur polynomials and power
sums by means of the characters of the Hecke algebras in the spirit of the classical Frobenius formula.
\end{abstract}

{\bf AMS Mathematics Subject Classification, 2010:} 05E05, 20C08


\section{Introduction}

Let $V$ be a vector space over the field $\mathbb{C}$, $\dim_\C V=N$. The action of the enveloping algebra $U(gl(N))$ on the space $V$ can be extended to the tensor products
$V^{\ot n}$, $n\ge 2$, via the coproduct $\Delta:U(gl(N))\to U(gl(N))^{\ot 2}$, defined on the generators $X_i\in gl(N)$, $1\le i\le N$ in the usual way:
$$
\Delta(X_i)\to X_i\ot 1+1\ot  X_i.
$$

Also, on the space $V^{\ot n}$ an action of the symmetric group $\mathcal{S}_n$  is naturally  defined via permutations of the factors in this tensor product.

The Schur-Weyl duality means that the actions of $U(gl(N))$ and $\mathcal{S}_n$ commute with each other. Moreover, the following module isomorphism takes place:
\be
V^{\ot n}\simeq \bigoplus_{\la\, \vdash n} \, V^\lambda\ot M^\lambda, \label{SW}
\ee
where the direct sum runs over all partitions $\la=(\la_1\geq \dots \geq \la_n)$ of the integer $n$, $V^\lambda$ and $M^\la$ are irreducible $U(gl(N))$ and $\mathcal{S}_n$ modules
respectively, labelled by the partition $\la$.

Now, consider a set of commutative  indeterminates $\{x_1,\dots,x_N\}$. With any partition $\la\vdash n$ one associates some families of symmetric polynomials\footnote{  Namely, the
polynomials invariant under any permutation of the elements  $x_i,\,\, 1\leq i\leq N$.} in these indeterminates. We need two of these families: the Schur polynomials (functions)
$s_\la(x_1,\dots, x_N)$ and the product of power sums $p_{\lambda}(x_1,\dots, x_N) = p_{\lambda_1}p_{\lambda_2}\dots p_{\lambda_n}$ where $p_k = \sum_i x_i^k$. The famous Frobenius formula reads:
\be
p_{\nu} (x_1,\dots, x_N)=\sum_{\la\,\vdash N} \chi^\la_\nu \,s_\la (x_1,\dots, x_N),
\label{Fro}
\ee
where $\chi^\la_\nu$ is the character of the symmetric group $\mathcal{S}_N$ in the representation $M^\lambda$ evaluated on an element of cyclic type  $\nu\vdash N$.

In \cite{J} M. Jimbo  established a similar duality between the quantum group $U_q(gl(N))$ and the Hecke algebras of $A_{n-1}$ type.  Guided by this result, A. Ram \cite{R} exhibited a
$q$-version of the Frobenius formula with the reference to the works  of  A.M. Vershik and S.V. Kerov as well as  of  R.C. King and B.G. Wybourne where  similar results were obtained.
In fact, the quantum group $U_q(gl(N))$ itself does not actually enter the corresponding constructions. As A. Ram observed, his ``presentation ... allows one to avoid the quantum group completely". More precisely, upon replacing in (\ref{Fro}) the character of the symmetric group by that of the Hecke algebra, evaluated on some special elements, A. Ram showed   that
the Hall-Littlewood polynomials appeared in the left hand side of (\ref{Fro}). Lately, his results have been generalized to super-versions of the Schur-Weyl duality and of the $q$-Frobenius formula (see \cite{M}).

In the current paper we present a new version of the Schur-Weyl duality and $q$-Frobenius formula, valid for the so-called Reflection Equation (RE) algebras. Each RE algebra $\LL(R)$
is associated with a skew-invertible Hecke type braiding $R:V^{\otimes 2}\to V^{\otimes 2}$ (see Section 2), we call these braidings  {\it  Hecke  symmetries}.  Well known examples of Hecke symmetries are those coming from the quantum group $U_q(sl(N))$  and their super-analogues. 

However, there exist big families of other Hecke symmetries, including  those which are not deformations either of the usual flips or of the super-flips (see \cite{G}).

In any RE algebra $\LL(R)$ one can define quantum  analogues  of various symmetric functions, in particular, the Schur functions $s_\lambda$ and the power sums $p_k$, which are
central in the corresponding RE algebra\footnote{In \cite{IOP} these objects are introduced in more general Quantum Matrix algebras, associated with pairs of braidings $(R,F)$.
However, in general, the corresponding symmetric functions are not central. Bellow, we restrict ourselves to the case $F=R$, which corresponds to the RE algebras.
Also, see Remark 14.}.

The quantum symmetric functions  those  we are dealing with  are  defined as  certain  homogeneous polynomials in generators of the RE algebra $\LL(R)$ (see Section 3). Our treatment
of these objects as symmetric functions is based on the remarkable similarity of their algebraic properties and those  of  the corresponding classical counterparts. In particular, the
quantum Schur  functions $s_\lambda$ are subject to  the Littlewood-Richardson rule with the classical coefficients for their products (see \cite{GPS2}). Also, there exist $q$-analogues
of the Newton identities connecting the power sums $p_k$ and the quantum elementary symmetric functions $s_{(1^n)}$. It should be emphasized that if a Hecke symmetry $R$ tends
to the usual flip $P$ as $q\to 1$, the algebra $\LL(R)$ tends to the algebra $\mathrm{Sym}(gl(N))$ and all quantum symmetric functions turn  into their classical  analogues  expressed via
generators of the latter algebra.

The classical symmetric functions can be also expressed in a more habitual  form via the eigenvalues of the generating matrix $L$ of the algebra  $\mathrm{Sym}(gl(N))$. The point is
that a similar presentation of the quantum Schur functions and power sums is also possible. The role of the eigenvalues of the generating matrix $L$ of the algebra $\LL(R)$ in this case
is  played by the roots of the Cayley-Hamilton polynomial for this matrix. Note that the coefficients of this polynomial are cental. Thus, the roots of the Cayley-Hamilton polynomial
are treated  as  elements of an algebraic extension of the algebra $\LL(R)$  and are assumed to be central in it.

If the initial Hecke symmetry $R$ is {\em even}, the quantum power sums are nothing but Hall-Littlewood polynomials (up to minor modifications, exhibited in the last section). If $R$ is of
general type (see the next section), the eigenvalues are split in two sets --- even eigenvalues and odd ones. Then all quantum symmetric functions in the corresponding RE algebra $\LL(R)$, being expressed via these two sets, become super-symmetric polynomials in the sense of \cite{S}.  Upon rewriting all (super-)symmetric polynomials via the eigenvalues of $L$, it is
possible to present the $q$-Frobenius formula in  a  form similar to that from \cite{R}, if $R$ is even, or similar to that from \cite{M}, if $R$ is of  general type.

The paper is organized as follows. In the next section we recall some facts about Hecke algebras of $A_{n-1}$ type and their representation theory. Also, we introduce the so-called $R$-matrix representations and the $R$-traces. In section 3 we introduce the RE algebras $\LL(R)$ and the quantum symmetric polynomials. Then we define a representation of the RE algebra on the tensor powers of the basic vector space $V$ and exhibit the quantum version of the Schur-Weyl duality. In the last section we prove the $q$-Frobenius formula connecting the quantum
symmetric polynomials defined in the RE algebra $\LL(R)$. Besides, we introduce the eigenvalues of the generating matrices $L$ and express all these polynomials via the eigenvalues.  Finally,  we illustrate the $q$-Frobenius formula by two explicit examples, corresponding to  the  Hecke symmetries $R$, coming from $U_q(sl(2))$ and $U_q(sl(1|1))$, respectively.

\section{$A_{n-1}$ type Hecke algebras and Hecke symmetries}

In this section we recall some facts on the theory of $A_{n-1}$ type Hecke algebras and its $R$-matrix representations. Most of these facts concerning the Hecke algebras and their representation theory can be found  in the review \cite{OP} (also, see the references therein).

Let $n\ge 2$ be a positive integer.
\begin{definition}
\label{def:1}
\rm
The $A_{n-1}$ type Hecke algebra $\H_n(q)$ is a unital associative algebra over the complex field $\mathbb{C}$ finitely generated by the set of {\it Artin generators} $\tau_i$,
$1\le i\le n-1$, which satisfy the following relations:
\begin{eqnarray}
\tau_i\, \tau_{i+1}\, \tau_i=\tau_{i+1}\, \tau_i\, \tau_{i+1},&\quad& 1\le i \le n-2,\nonumber\\
\rule{0pt}{4mm}
\tau_i\, \tau_j=\tau_j\, \tau_i\, &\quad& |i-j|\geq 2,\label{H-alg}\\
\rule{0pt}{4mm}
(\tau_i-q\, e)(\tau_i+q^{-1} e)=0,&\quad& 1\le  i\le n-1,\nonumber
\end{eqnarray}
where $e$ is a unit element of the algebra and $q\in\mathbb{C}\setminus \{0,\pm 1\}$ is a numeric parameter.
\end{definition}

Note that at $q=\pm 1$ the defining relations on the generators $\tau_i$ turn into those on the generators of the symmetric group $\S_n$. For $q\not=\pm 1$ the Hecke algebra $\H_n(q)$
is a deformation of the group algebra $\C[\S_n]$.

The following proposition establishes some important properties of the Hecke algebra $\H_n(q)$.

\begin{proposition}
\label{prop:2}
\phantom{a}
\begin{enumerate}
\item The Hecke algebra $\H_n(q)$ is finite dimensional: $\mathrm{dim}_{\,\C} \H_n(q) = n!$.
\item If for any integer $2\le k\le n$ the parameter $q$ satisfies the condition $q^{2k}\not=1$,  then the algebra $\H_n(q)$ is semisimple and isomorphic to the group algebra of the
symmetric group $\S_n$:
\be
\H_n(q)\simeq \C[\S_n].
\label{HS-iso}
\ee
\end{enumerate}
\end{proposition}

The values of $q$ which do not belong to the discrete set $q^{2k} =1$, $2\le k\le n$, (where $\H_n(q)$ loses the semisimplicity) are called {\it generic}. For generic values of $q$ the
$q$-analogues of integers do not vanish:
$$
k_q:=\frac{q^k-q^{-k}}{q-q^{-1}} \not=0,\quad 1\le \forall\,k\le n.
$$
In what follows we shall always assume the parameter $q$ to be generic.

Now, we turn to the consequence of the semisimplicity of the algebra $\H_n(q)$. According to the Wedderburn-Artin theorem, the finite dimensional semisimple Hecke algebra
$\H_n(q)$ is isomorphic to a direct product of matrix algebras over $\C$. Besides, in virtue of the isomorphism (\ref{HS-iso}) the matrix factors can be labelled by partitions of $n$:
\be
\H_n(q)\simeq M_{\lambda(1)}(\C)\times M_{\lambda(2)}(\C)\times \dots\times M_{\lambda(s)}(\C),
\label{H-dir-pr}
\ee
where $\{\lambda(1),\lambda(2),\dots ,\lambda(s)\}$ is a (unordered) set of {\it all} possible partitions of the integer $n$: $\lambda(k)\vdash n$. The dimension of the
component $M_\lambda(\C)$ equals  $(d_\lambda)^2$, where $d_\lambda$ is the number of the {\it standard Young tables} corresponding to the Young diagram of the partition $\lambda$.

In each matrix algebra $M_{\la}$ there exists a basis of matrix units\footnote{Recall, that a matrix unit $E_{ij}$ is a matrix with the only nonzero element 1 located at the crossing of
the $i$-th row and $j$-th column.} $\{E^\la_{ij}\}_{1\le i,j\le d_\lambda}$. Denote $e^\la_{ij}$ the element of the Hecke algebra $\H_n(q)$, which corresponds to $E^\la_{ij}$ under the isomorphism (\ref{H-dir-pr}). Note that the multiplication table for these elements reads:
\be
e^\la_{ki}\,e^\mu_{r p}=\delta^{\la\mu}\delta_{ir}\,e^\la_{kp}.
\label{prod-mun}
\ee
Evidently, the full set of elements $\{e^\lambda_{ij}\}_{1\le i,j\le d_\lambda}$ for all $\lambda\vdash n$ forms a linear basis of the Hecke algebra $\H_n(q)$.

Now, consider the left regular representation $\pi$ of the Hecke algebra $\H_n(q)$ in the basis $\{e^\la_{ki}\}$. Let $z\in \H_n(q)$ be an arbitrary element. Its image $\pi(z)$ in the basis
set $\{e^\la_{ki}\}$ for each fixed $\lambda\vdash n$ is represented by a $d_\la\times d_\la$ matrix $Z^\lambda$:
\be
\pi(z) e^\la_{ki}\stackrel{\mbox{\rm\tiny def}}{=}z\,e^\la_{ki}=\sum_{r=1}^{d_\lambda} Z_{kr}^\la\,  e^\la_{ri}.
\label{maa}
\ee

The elements $e^\la_{ii}$, $1\le i \le d_\la$, are primitive orthogonal idempotents defining a resolution of the unit $e$:
\be
e=\sum_{\lambda\,\vdash n}\sum_{i=1}^{d_\lambda}e^{\lambda}_{ii}, \qquad e^{\lambda}_{ii}\,e^{\mu}_{jj} = \delta^{\lambda\mu}\delta_{ij}\,e^{\lambda}_{ii}.
\label{reso}
\ee

In \cite{OP} there is exhibited an explicit construction of the elements $e^\la_{ij}$ via the generators of the Hecke algebra $\H_n(q)$. Namely, each idempotent $e^\lambda_{ii}$
is a polynomial in the so-called {\it Jucys-Murphy elements} $\{j_k\}_{1\le k\le n}$, which are defined by recursion as follows:
\be
j_1 = e,\qquad j_k = \tau_{k-1}\,j_{k-1}\tau_{k-1},\quad 2\le k\le n.
\label{JM}
\ee
Note that the set of Jucys-Murphy elements generates a maximal commutative subalgebra in $\H_n(q)$. Besides, the center of $\H_n(q)$ consists of all symmetric polynomials in the
Jucys-Murphy elements.

\begin{remark} \rm The idempotents $e^\la_{ii}$ constructed in \cite{OP} are common eigenvectors of the commutative set $\{j_k\}_{1\le k\le n}$ with respect to the left or right
multiplication:
$$
j_k\,e_{ii}^{\lambda} = e_{ii}^{\lambda}\, j_k = q^{2c_i(k)}e_{ii}^{\lambda}.
$$
Here $c_i(k)=x_i(k)-y_i(k)$ is the {\it content} of a box in the standard Young table $(\lambda,i)$ in which the integer $k$ is placed, $x_i(k)$ (respectively $y_i(k))$ is the number
of the column (the row) of the box where $k$ is located.
\end{remark}

Our next objective is to introduce  the so-called $R$-matrix representations of the Hecke algebras.

Given a finite dimensional complex vector space $V$,  $\mathrm{dim}_{\,\C}\, V= N$, we call an operator $\hat R\in \mathrm{End}(V^{\otimes 2})$  {\it a braiding}
if it satisfies the following equation in the algebra $\mathrm{End}(V^{\otimes 3})$:
$$
(\hat R\otimes\mathrm{Id}_V)(\mathrm{Id}_V\otimes \hat R)(\hat R\otimes\mathrm{Id}_V) = (\mathrm{Id}_V\otimes \hat R)(\hat R\otimes\mathrm{Id}_V) (\mathrm{Id}_V\otimes \hat R).
$$
If we fix a basis $\{x_i\}_{1\le i\le N}$ in the space $V$, then in the corresponding tensor basis $\{x_i\otimes x_j\}$ of $V^{\otimes 2}$ the braiding $\hat R$ is described by an
$N^2\times N^2$ matrix $R = \|R_{i j}^{k l}\|$ and the above equation on $\hat R$ can be rewritten in the matrix form:
\be
R_{12}R_{23}R_{12}  = R_{23}R_{12}R_{23},
\label{YB}
\ee
where $R_{12} = R\otimes I$, $R_{23} = I\otimes R$, $I$ is the $N\times N$ unit matrix. The equation (\ref{YB}) is known as the {\it braid relation}.

Below, we use the following notation
\be
R^{\,\pm 1}_{\,i} = I^{\otimes (i-1)}\otimes R^{\,\pm 1}\otimes I^{\otimes (m-i-1)},\quad 1\le i\le m-1.
\label{embed}
\ee
Therefore, each $R^{\,\pm 1}_{\,i}$ defines an embedding of $R^{\,\pm 1}$ into the space of $N^m\times N^m$ matrices for an arbitrary given integer $m\ge 2$.

Among the matrix solutions $R$ of the braid relation (\ref{YB}) we are interested in those satisfying an additional {\it Hecke condition}:
\be
R^{2} = I\otimes I +(q-q^{-1})R,\qquad q\in \C\setminus 0.
\label{Hec}
\ee
Such  solutions will be referred to as {\it Hecke symmetries}.

Note that at $q=1$ the relation (\ref{Hec}) gives $R^2=I\ot I$. The braidings meeting this relation will be called {\it involutive symmetries}. As examples of involutive symmetries we can
mention the usual flip $P$ and the super-flips $P_{m|n}$, where $m$ and $n$ are respectively the dimensions of the even $V_0$ and odd $V_1$ subspaces of the super-space
$V=V_0\oplus V_1$. By the super-dimension of this space we mean the ordered pair $(m|n)$.

As a generalization of this notion we introduce a bi-rank of an involutive or Hecke symmetry $R$. For this purpose we consider an $R$-skew-symmetric algebra of the space $V$
(see \cite{GPS3} for detail) and the corresponding  Hilbert-Poincar\'e series. According to \cite{H} this series is always a ratio of two coprime polynomials. The ordered pair of integers
$(r|s)$ where $r$ and $s$ are respectively the degrees of  polynomials in the numerator and denominator of the Hilbert-Poincar\'e series is called {\it the bi-rank} of the symmetry $R$.
If $s=0$, the corresponding symmetry $R$ is called {\em even}.

Besides the Hecke condition (\ref{Hec}) we shall always assume the symmetries $R$ to be skew-invertible.

\begin{definition}\rm
\label{def:4}
Let  $R$ be an involutive or a Hecke symmetry. We say that it is {\it a skew-invertible} symmetry if there exists an operator $\Psi:V^{\otimes 2}\to V^{\otimes 2}$ such that
$$
\Tr_2 R_{12}\,  \Psi_{23}=\Tr_2 (R\ot I)\, (I\ot \Psi)=P_{13}\quad \Leftrightarrow \quad \sum_{r,j\,=1}^NR_{ij}^{kr}\,\Psi_{rm}^{jn}=P_{im}^{kn}=\delta_m^k \delta_i^n.
$$
\end{definition}

Hereafter, we identify the operators $\Psi$ and the corresponding matrix in the fixed basis.

If $R$ is a skew-invertible involutive or Hecke symmetry, we introduce the $N\times N$ matrices $B=\|B_i^j\|$ and $C = \|C_i^j\|$ as partial traces of the matrix $\Psi$:
\be
B_i^j = \sum_{a=1}^N \Psi_{ai}^{aj},\qquad C_i^j = \sum_{a=1}^N \Psi_{ia}^{ja}.
\label{B-C}
\ee
If the bi-rank of the symmetry $R$ is $(r|s)$ then the following relations hold \cite{GPS3}:
\be
BC= q^{-2(r-s)}\, I,\qquad \Tr\, B=\Tr\, C= q^{s-r}\,(r-s)_q.
\label{form}
\ee
The matrix $B$ is used for construction of representations of the RE algebra $\LL(R)$ (see the next Section) while the matrix $C$ enter the definition of the $R$-trace.
\begin{definition}\rm
Let $\mathrm{Mat}_N(U)$ be the space of $N\times N$ matrices with entries from an arbitrary complex vector space $U$. The quantum or $R$-trace is the map
$\Tr_R: \mathrm{Mat}_N(U)\rightarrow U$ defined by the rule:
\be
\Tr_R\, M:=\Tr\,( C\cdot M)\qquad \forall \, M\in \mathrm{Mat}_N(U).
\label{tra}
\ee
\end{definition}
Note that the $R$-trace plays an important role in studying the structure of the RE algebra $\LL(R)$, in particular, in definitions of quantum symmetric polynomials related to the RE
algebras.

Below we use also the multiple $R$-trace defined on the tensor product of matrix spaces. Let $M_{12\dots k}$ be an $N^k\times N^k$ matrix from the space
$\mathrm{Mat}_N(U)^{\otimes k}$. Then the multiple $R$-trace is defined by the rule:
$$
\Tr_{R(12\dots k)}M_{12\dots k} = \Tr_{(12\dots k)}(C_1C_2\dots C_kM_{12\dots k})
$$
where  $C_i = I^{\otimes(i-1)}\otimes C\otimes I^{\otimes(k-i)}$.

At last, we present a few technical formulae needed for further  calculations. As a straightforward consequence of Definition \ref{def:4} and (\ref{B-C})
we have
$$
\Tr_{R(2)}R_1 = I\quad \Leftrightarrow \quad \sum_{a,b\,=1}^NC_a^bR_{ib}^{ja} = \delta_i^j.
$$
A direct generalization of this relation reads:
\be
\Tr_{R(k+1)}\, R_k=I^{\otimes k}.
\label{trrr}
\ee
Besides, it immediately follows from (\ref{form}) that:
$$
\Tr_R I = \Tr \,C = q^{s-r}(r-s)_q.
$$

Now, we introduce a special representation of the Hecke algebra. Given a Hecke symmetry $R$, we define the {\it $R$-matrix representation} $\rho_R$ of the Hecke algebra
$\H_n(q)$ in the space $V^{\otimes n}$ by sending the Artin generator $\tau_i$ to the matrix $R_{\,i}$ (\ref{embed}):
\be
\rho_R(\tau_i) = R_{\,i} = I^{\otimes (i-1)}\otimes R\otimes I^{\otimes (n-i-1)},\qquad 1\le i\le n-1.
\label{R-rep}
\ee
This representation plays a key role in the definition of the so-called characteristic subalgebra of the RE algebra which, in particular, contains the quantum symmetric functions.

\section{Schur-Weyl duality and symmetric polynomials}

From now on we assume $R:V^{\ot 2}\to V^{\ot 2}$, $\dim_\C\, V=N$, to be a skew-invertible Hecke symmetry of  bi-rank $(r|s)$.

\begin{definition}\rm
The Reflection Equation (RE) algebra $\LL(R)$ is a unital associative algebra over the complex field $\C$ generated by $N^2$ entries of the matrix $L=\|l_i^j\|_{1\leq i, j \leq N}$
subject to the relation:
\be
R\, L_1 R\, L_1-L_1 R\,L_1 R=0,\qquad L_1=L\ot I.
\label{RE}
\ee
The matrix $L$ is called the generating matrix.
\end{definition}

Our next objective is to describe a finite dimensional representations of the RE algebra in the spaces $V^{\ot n}$. Note that if $R$ comes from the quantum group $U_q(sl(N))$ and $q$
is generic, then the category of finite dimensional representations of the RE algebra is just the same as that of $U_q(sl(N))$. However, the coalgebraic structures of $U_q(sl(N))$ and
RE algebra $\LL(R)$ are completely different: while the quantum group has the usual bi-algebraic structure, in the RE algebra it is in a sense braided (see \cite{GPS3} for detail). Our method
of constructing the mentioned category is universal and valid for all RE algebras, associated with any skew-invertible Hecke symmetries.

Let us define a linear action of the generators $l_i^j$ of the RE algebra on a basis vector $x_k$ of the space $V$ by the following rule:
\be
l_i^j \triangleright x_{ k} =\delta_i^j\, x_{ k}-(q-q^{-1})  B_k^j\, x_{ i}.
\label{actt}
\ee
Hereafter, $\triangleright$ denotes the action of a linear operator on a vector, the matrix $B=\|B_i^j\|$ was introduced in (\ref{B-C}). It is not difficult to show that this action defines a representation of the algebra $\LL(R)$ in the space $V$. The action (\ref{actt}) can be written in the following matrix form:
$$
L_2\triangleright x_{|1\rangle}=(I\ot I)\, x_{|1\rangle}-(q-q^{-1}) B_1 P_{12} x_{|1\rangle},
$$
where $L_2 = I\otimes L$. In order to make this matrix writing more transparent, we present the last term of this formula with the use of the indexes, numerating  elements of the
corresponding vector or matrix spaces. Here we use the ``matrix summation rule'': if a relation contains the product of two objects with the same number of matrix space then it is
assumed a matrix product in this space. For example, the explicit matrix indices of the second term read:
$$
\left(B_1 P_{12}\, x_{|1\rangle}\right)_{i_1i_2}^{\,\,\,j_2} = \sum_{a,b\,=1}^NB_{i_1}^aP_{ai_2}^{\, bj_2}\,x_b.
$$

For the skew-invertible Hecke symmetry one can rewrite the above relation in an {\it equivalent} form:
$$
L_1R_1\triangleright x_{|1\rangle}= R^{-1}_1  x_{|1\rangle}.
$$

Now we introduce the embeddings of the generating matrix $L\in \mathrm{Mat}_N(\LL(R))$ into the space $\mathrm{Mat}_N(\LL(R))^{\otimes m}$:
\be
L_{\ov 1}=L_{\un 1} = L\otimes I^{\otimes (m-1)},\qquad L_{\ov {k+1}}=R_k L_{\ov k}R_k^{-1},\qquad L_{\un {k+1}}=R_k^{-1}L_{\un k}R_k,\quad 1\le k\le m-1,
\label{over-not}
\ee
where $R^{\pm 1}_i$ is defined in (\ref{embed}). The concret value of the integer $m$ does not matter: it should be just sufficiently large in order all the matrix relations would make sense.

With notation (\ref{over-not}) the above matrix form of the RE algebra action can be further transformed:
\be
L_1R_1\triangleright x_{|1\rangle}= R^{-1}_1  x_{|1\rangle}\quad\Leftrightarrow\quad
R^{-1}_1 L_1R_1 \triangleright x_{|1\rangle}=R^{-2}_1 x_{|1\rangle}\quad\Leftrightarrow \quad L_{\un 2} \triangleright x_{|1\rangle} =J_2^{-1} x_{|1\rangle},
\label{obt}
\ee
where $J_2 = R_1^2 = \rho_R(j_2)$ is the image of the Jucys-Murphy element $j_2$ (\ref{JM}) under the $R$-matrix representation (\ref{R-rep}). The images
of other Jucys-Murphy elements will be denoted by analogous symbols: $J_k=\rho_R(j_k)$.

The representation of the RE algebra in the tensor powers $V^{\otimes k}$ is given by a generalization of the formula (\ref{obt}). The corresponding result was obtained in \cite{GPS5},
we reproduce it here for the reader's convenience.

\begin{proposition}{\rm \bf \cite{GPS5}}
The representation of the RE algebra $\LL(R)$ in the space $V^{\otimes k}$ is defined by the following action of its generators on the tensor basis of the space $V^{\otimes k}$:
\be
L_{\un{k+1}}\triangleright (x_{|1\rangle}\ot\dots \ot x_{|k\rangle})=J_{k+1}^{-1}\,x_{|1\rangle}\ot\dots \ot x_{|k\rangle},
\label{tri}
\ee
where $J_{k+1} = \rho_R(j_{k+1}) = R_k\dots R_2R_1^2R_2\dots R_k$.
\end{proposition}

\noindent
{\bf Proof.}  First, observe that the defining relations (\ref{RE}) of the RE algebra $\LL(R)$ can be written in the equivalent form:
$$
R_1L_{\un 2}L_{\un 1}=L_{\un 2}L_{\un 1}R_1.
$$
Moreover, they can be presented in terms of higher copies of $L$ defined in (\ref{over-not}):
\be
R_{k+1} L_{\un{k+2}}L_{\un{k+1}}=L_{\un{k+2}}L_{\un{k+1}} R_{k+1}, \qquad \forall\, k\ge 0.
\label{shi}
\ee

Now we apply the left hand side of (\ref{shi}) to the tensor basis $x_{|1\rangle}\ot\dots \ot x_{|k\rangle}$:
\begin{eqnarray*}
R_{k+1}\hspace*{-6mm}&&\hspace*{-4mm}L_{\un{k+2}}L_{\un{k+1}}\triangleright x_{|1\rangle}\ot\dots \ot x_{|k\rangle}=R_{k+1} L_{\un{k+2}}\triangleright (J_{k+1}^{-1}
 x_{|1\rangle}\ot\dots \ot x_{|k\rangle})\\
&=&\hspace*{-3mm}R_{k+1}J_{k+1}^{-1}(R_{k+1}^{-1} L_{\un{k+1}} R_{k+1})\triangleright  x_{|1\rangle}\ot\dots \ot x_{|k\rangle}=R_{k+1}J_{k+1}^{-1}R_{k+1}^{-1} J_{{k+1}}^{-1}R_{k+1}
x_{|1\rangle}\ot\dots \ot x_{|k\rangle}\\
&=&\hspace*{-3mm}R_{k+1}J_{k+1}^{-1}J_{k+2}^{-1}R_{k+1}^2 x_{|1\rangle}\ot\dots \ot x_{|k\rangle}.
\end{eqnarray*}
Here, in the second line of calculations we used the fact that the matrices $L_{\un{k+2}}$ and $J_{k+1}^{-1}$ commute with each other since the $J^{-1}_{k+1}$ is a polynomial in matrices
$R^{-1}_{\,i}$, $1\le i\le k$. This commutativity is a simple consequence of the braid relation on $R$ and definition (\ref{over-not}) of $L_{\underline{k+2}}$.

Next, we apply  the right hand side of (\ref{shi}) to the same basis and get:
\begin{eqnarray*}
L_{\un{k+2}}\hspace*{-6mm}&&\hspace*{-4mm}L_{\un{k+1}}R_{k+1}\triangleright x_{|1\rangle}\ot\dots \ot x_{|k\rangle}=J_{k+1}^{-1}L_{\underline{k+2}}\triangleright
x_{|1\rangle}\ot\dots \ot x_{|k\rangle} R_{k+1}\\
&=&\hspace*{-3mm}J_{k+1}^{-1}R_{k+1}^{-1}J_{k+1}^{-1} R_{k+1}^2 x_{|1\rangle}\ot\dots \ot x_{|k\rangle} = R_{k+1} J_{k+2}^{-1}J_{k+1}^{-1} R_{k+1}^2
x_{|1\rangle}\ot\dots \ot x_{|k\rangle}.
\end{eqnarray*}
Here we used the fact that the matrix $R_{k+1}$ does not affect the product $x_{|1\rangle}\ot\dots \ot x_{|k\rangle}$.

The actions of both sides of (\ref{shi}) on basis vectors of $V^{\otimes k}$ differ only by the order of matrices $J_{k+1}^{-1}$ and $J_{k+2}^{-1}$ which commute with each
other and, therefore, these actions are identical. So, the action (\ref{tri}) gives a representation of the RE algebra in $V^{\otimes k}$. \hfill \rule{6.5pt}{6.5pt}

\medskip

The Schur-Weyl duality in our setting is formulated in the following proposition.

\begin{proposition}
The action of the Hecke algebra $\H_k(q)$ onto the space $V^{\ot k}$, defined via the representation $\rho_R$ {\rm (\ref{R-rep})}, commutes with the action {\rm (\ref{tri})} of the RE
algebra $\LL(R)$.
\end{proposition}

\noindent
{\bf Proof.} The claim follows immediately from the fact that the Jucys-Murphy element $j_{k+1}$ commutes with generators $\tau_i$, $1\le i\le k-1$, and
consequently the operator  $L_{\underline{k+1}}$ from (\ref{tri})
commutes with  $R_i$, $1\le i\le k-1$.\hfill\rule{6.5pt}{6.5pt}

\medskip

Our next objective is to define some quantum symmetric polynomials on the algebras $\LL(R)$. Let us introduce the shorthand notation
$$
L_{\ov{1\to k}}=L_{\ov 1}L_{\ov 2}\dots L_{\ov k}.
$$
Due to defining relations of the RE algebra (\ref{shi}) we have the equality:
\be
R_{k-1}R_{k-1}\dots R_1L_{\ov{1\to k}} = L_{\ov{1\to k}}\,R_{k-1}R_{k-1}\dots R_1.
\label{gen-shi}
\ee
\begin{proposition} {\rm\bf \cite{IP}}
Let $z\in \H_n(q)$ be an arbitrary element. Then the element
\be
ch_n(z)=\Tr_{R(1\dots n)}\left(\rho_R(z) L_{\ov{1\to n}}\,\right) = \Tr_{R(1\dots n)} \left( L_{\ov{1\to n}}\,\,\rho_R(z)\right)
\label{sem}
\ee
is central in the algebra $\LL(R)$.
\end{proposition}

Note that the second equality in {\rm (\ref{sem})} is fulfilled in virtue of the cyclic property\footnote{This property means that
$\Tr_{R(1\dots k)} \,A\, B=\Tr_{R(1\dots k)}\, B\, A$ where $A$ and $B$ are two $N^k\times N^k$ matrices, provided one of them is a polynomial in $R_1,\dots, R_{k-1}$.
This property is due to the relation $R_iC_iC_{i+1}=C_iC_{i+1} R_i$ (see \cite{GPS3}).} of the $R$-trace or  formula (\ref{gen-shi}).

Thus, we have defined a map
$$
ch_n: \H_n(q)\to Z(\LL(R)),\quad z\mapsto ch_n(z),
$$
which is  called {\em characteristic}.
Here, $Z(A)$ stands for the center of the algebra $A$.

\begin{remark}\rm
In the work \cite{IP} the subalgebra generated by all the elements $ch_n(z)$ for all $n\in \mathbb{Z}_+$ was called the {\it characteristic subalgebra} and it was defined for a general
quantum matrix algebra associated with a pair of compatible braidings $(R,F)$. In the general case the characteristic subalgebra was proved to be commutative but not central.
\end{remark}

Now,  consider two particular examples of the elements $z\in \H_k(q)$, giving rise to the quantum symmetric polynomials. First, we set $z_{(k)}=\tau_{k-1}\tau_{k-2}\dots \tau_1$
($z_{(k)}$ is called the Coxeter element). Then we get
\be
p_k(L):=\Tr_{R(1\dots k)}\left(\rho_R(\tau_{k-1}\tau_{k-2}\dots \tau_1) L_{\ov{1\to k}}\,\right)=\Tr_{R(1\dots k)}\left( R_{k-1} R_{k-2}\dots R_1L_{\ov{1\to k}}\,\right).
\label{pws}
\ee
These elements are called {\em power  sums}.

Note that these power sums can be reduced  to $p_k(L)=\Tr_R\, L^k$. Indeed, with the use of the cyclic property of the $R$-trace or, alternatively, the algebraic
relation (\ref{gen-shi}) we rewrite $p_k = \Tr_{R(1\dots k)}(L_{\ov{1\to k}}\,\, R_{k-1} R_{k-2}\dots R_1)$ and then use the definition of $L_{\overline k}$ and the formula (\ref{trrr}).

For any partition $\la=(\la_1,\dots,\la_k)$ we introduce the corresponding symmetric function  $p_\lambda(L)$ in the standard way, that is as the following product:
\be
p_{\la}(L) = p_{\la_1}(L)\dots p_{\la_k}(L).
\label{gen-ps}
\ee
Since all factors in this product are central in the algebra $\LL(R)$ their order does not matter.
We  call the symmetric function $p_\lambda(L)$ power sum, corresponding to the partition $\la$. Thus,
the power sums $p_k(L)$ is a particular case corresponding to the one-row diagrams $\la=(k)$.

In the second example we define  the {\it quantum Schur polynomials} associated with any partition $\lambda\vdash k$, $k\ge 1$ as follows:
\be
s_\la(L)=\Tr_{R(1\dots k)}\left( e^\la_{ii}(R)\,L_{\ov{1\to k}}\,\right).
\label{schu}
\ee
Hereafter, we use the notation $e^\la_{ij}(R)=\rho_R(e^\la_{ij})$.

The quantum Schur polynomials were introduced in \cite{GPS1}, where it was shown that the right hand side of (\ref{schu}) depends only on the partition $\lambda$ but
not on the number $i$ of the Young table used in construction of the idempotent $e_{ii}^\lambda$.

Above, we have already mentioned the similarity of the quantum symmetric polynomials and their classical analogues. Observe once more that for the quantum Schur polynomials  the
Littlewood-Richardson rule
$$
s_\lambda(L)s_\mu(L) = \sum_{\nu}C_{\lambda\mu}^\mu s_\nu(L),
$$
is valid. Here  positive integers $C_{\lambda\mu}^\mu$ are the {\it classical} Littlewood-Richardson coefficients. Besides, the quantum Schur polynomials form a linear basis in the central characteristic subalgebra of the RE algebra $\LL(R)$.

Turn now to the structure of the elements (\ref{pws}). Each $p_k(L)$ is defined with the use of the Coxeter element $ z_{(k)}=\tau_{k-1}\, \tau_{k-2}\dots \tau_1\in \H_k(q)$. We claim
that any product (\ref{gen-ps}) of the power sums parameterized by a partition $\nu\vdash k$ can be presented as the $R$-trace of the form (\ref{pws}) with special element
$z_\nu\in \H_k(q)$:
$$
p_{\nu} = \Tr_{R(1\dots k)}\left(\rho_R(z_\nu)L_{\overline{1\rightarrow k}}\,\right).
$$

For this purpose we introduce the {\it Coxeter element with gaps} $z_\nu$ of the cyclic type $\nu\vdash k$. We start from the Coxeter element $z_{(k)}$ and remove some
$s\le k-1$ factors $\tau_i$ from the string $\tau_{k-1}\tau_{k-2}\dots \tau_1$. The element $z_\nu\in \H_n(q)$ obtained in this way is the product of a few strings of  generators
like $\tau_{j_1}\,\tau_{j_1-1}\dots \tau_{j_2}$ with some $j_2\le j_1$. The number $\ell_{j_1} = j_1-j_2+1$ is called the length of the string. Note, that the strings commute with each
other since the difference of the indices of generators constituting two arbitrary strings is greater than 1. This element $z_\nu$ will be called the Coxeter element with gaps.

The {\it cyclic type} of $z_\nu$ is a partition $\nu\vdash k$ defined as follows. We sort the string lengths of $z$ in non-increasing order: $\ell_{i_1}\ge \ell_{i_2}\ge\dots\ge \ell_{i_r}\ge 1$ and
then we set by definition:
$$
\nu_1 = \ell_{i_1}+1,\quad \nu_2 = \ell_{i_2}+1,\quad\dots\quad  \nu_r = \ell_{i_r}+1.
$$
If the sum of components $\nu_i$ turns out to be less than $k$
$$
k-\sum_{i=1}^r\nu_i = p\ge 1,
$$
then we must complete the set of components of the partition $\nu$ by $p$ additional components equal to 1: $\nu_{r+1} = \dots=\nu_{r+p} = 1$.

\medskip

\noindent{\bf Example.} Let $z_{(7)}= \tau_6\tau_5\tau_4\tau_3\tau_2\tau_1\in \H_7(q)$. Cancelling $\tau_5$ and $\tau_1$ we get the element $z_\nu = \tau_6\tau_4\tau_3\tau_2$.
The element $z_\nu$ consists of two strings: $\tau_4\tau_3\tau_2$ of the length 3 and $\tau_6$ of the length 1. So, we find two components $\nu_1 =4$ and $\nu_2 = 2$ of the
cyclic type $\nu\vdash 7$. But since $\nu_1+\nu_2 = 6<7$ we have to add one more component $\nu_3 =1$. Therefore, the cyclic type of $z_\nu = \tau_6\tau_4\tau_3\tau_2$
is the partition $\nu =(4,2,1)$. If we cancel the generator $\tau _4$, we get $z_\mu = \tau_6\tau_5\tau_3\tau_2\tau_1$ of the cyclic type $\mu=(4,3)$.

From the other hand, it is easy to see that inverse construction is always possible: given any partition $\lambda\vdash k$ we can construct a Coxeter element with gaps $z_\lambda$
with the cyclic type $\lambda$ starting from $z_{(k)}$. Note, that this construction is not unique, but all the elements $z_\lambda$ of the same cyclic type are conjugated in the Hecke
algebra $\H_k(q)$ (even in the braid group $\mathcal{B}_k$) and, therefore, leads to the same power sum $p_{\lambda}(L)$.

\begin{proposition}
\label{prop:11}
Let $z_\nu \in \H_k(q)$ be a Coxeter element with gaps of the cyclic type $\nu=(\nu_1,\dots ,\nu_k)$. Then the following holds
$$
ch_k(z_\nu)=\Tr_R( L^{\nu_1})\dots \Tr_R (L^{\nu_k}) = p_{\nu}(L).
$$
\end{proposition}

\noindent
{\bf Proof.} According to definition (\ref{sem}):
$$
ch_k(z) =\Tr_{R(1\dots k)}( L_{\ov{1\to k}}\,  \rho_R(z_\nu)).
$$
If the factor $\tau_{k-1}$ is removed in the element $z_\nu$ at the position $k$ we have
$$
\Tr_{R(k)}(L_{\ov{k}})= \Tr_{R(1)} (L_{\ov{1}})
$$
in virtue of the following property of the $R$-trace
$$
\Tr_{R(k)} L_{\ov k}= I^{\otimes (m-1)} \Tr_{R}(L),\quad 1\le \forall k\le m,
$$
where in accordance with definition (\ref{over-not}) $L_{\overline k}$ is treated as the matrix from $\mathrm{Mat}_N(\LL(R))^{\otimes m}$ for some sufficienly large $m$.

If $\tau_{k-1}$ is not removed, we use the relation
$$
L_{\ov k}\, R_{k-1}=R_{k-1} L_{\ov{k-1}}
$$
and apply (\ref{trrr}). Thus, we get
$$
ch_k(z_\nu) =
\Tr_{R(1\dots k-1)}\left( L_{\ov{1\to k-1}}\,L_{\ov{k-1}}\, \rho_R(z_{\mu})\right)=\Tr_{R(1\dots k-1)}\left(L_{\ov{1\to k-2}}\,L_{\ov{k-1}}^2\, \rho_R(z_{\mu})\right),
$$
where $z_\mu$ is obtained from $z_\nu$ by removing the factor $\tau_{k-1}$.

By continuing this procedure, we extract the multilier $\Tr_R (L^{\nu_1})$ in the structure of $ch_k(z_\nu)$, if the string, beginning by the factor $\tau_{k-1}$,  has the length $\nu_1-1$.
The same consideration is applicable for any string of the length $\ell_i =\nu_i-1$: it gives rise to the multiplier $\Tr_R(L^{\nu_i})$.

At last, each component of $\nu$ which is equal to 1 (if it exists) corresponds to the absence of some generator $\tau$ and gives rise to the multiplier $\Tr_R(L)$.
\hfill\rule{6.5pt}{6.5pt}

\section{Two forms of the $q$-Frobenius formula}
\label{sec:4}

First, we exhibit the $q$-Frobenius formula for the quantum symmetric polynomials realized in terms of generators of the RE algebra $\LL(R)$.

\begin{proposition}
Let $\nu$ be a partition of a given positive integer $n$. Then the following relation holds:
\be
p_{\nu}(L)= \sum_{\la\,\vdash n}\chi_\nu^\la \,s_\la(L),
\label{q-Fr}
\ee
where $\chi_\nu^\la$ is the character of the Hecke algebra $\H_n(q)$ in the representation $M^\lambda$ evaluated on an element of cyclic type $\nu$.
\end{proposition}

\noindent
{\bf Proof.}
Let $z$ be an arbitrary element of the Hecke algebra  $\H_n(q)$.  Then we have
\be
\Tr_{R(1\dots n)}\left(\rho_R(z)\,L_{\ov{1\to n}}\,\right) =\Tr_{R(1\dots n)}\Big( \rho_R(z) \sum_{\lambda\,\vdash n}\sum_{i=1}^{d_\lambda}
e^{\lambda}_{ii}(R)\,L_{\ov{1\to n}}\Big),
\label{exp}
\ee
where the resolution of the unity (\ref{reso}) is used. Now, we apply the left regular representation (\ref{maa}) of $\H_n(q)$. Thus, we get:
$$
\rho_R(z)e_{ii}^\lambda(R) = \rho_R(z \,e_{ii}^\lambda) = \sum_{k=1}^{d_\lambda}Z_{ik}^\lambda\, e_{ki}^\lambda(R).
$$
Now, we can  rewrite (\ref{exp}) as follows
$$
\Tr_{R(1\dots n)}\left(\rho_R(z)\,L_{\ov{1\to n}}\,\right) = \sum_{\lambda\,\vdash n}\sum_{i,k=1}^{d_\lambda}Z^\lambda_{ik}\,\Tr_{R(1\dots n)}(e_{ki}^\lambda(R)\,
L_{\overline{1\rightarrow n}}\,).
$$

Next, we use the identity:
$$
\Tr_{R(1\dots n)}(e_{ki}^\lambda(R)\, L_{\overline{1\rightarrow n}}\,) = \delta_{ik}\,s_\lambda(L).
$$
Indeed, as follows from the multiplication table of matrix units (\ref{prod-mun}): $e^\lambda_{ki}(R) = e^\lambda_{ki}(R)e^\lambda_{ii}(R)$. Besides, the matrices $e^\lambda_{ii}(R)$
are polynomials in $R_i$, $1\le i\le n-1$, and, therefore, they commute with $L_{\overline{1\rightarrow n}}$ (see (\ref{gen-shi})) and can be moved under the $R$-trace by cyclic
property. All this facts allows us to make the following identical transformations:
\begin{eqnarray*}
\Tr_{R(1\dots n)}(e_{ki}^\lambda(R) L_{\overline{1\rightarrow n}}\,) \hspace*{-2.5mm}&=&\hspace*{-2.5mm} \Tr_{R(1\dots n)}(e_{ki}^\lambda(R)
L_{\overline{1\rightarrow n}}\,e_{ii}^\lambda(R)) =  \Tr_{R(1\dots n)}(e_{ii}^\lambda(R)e_{ki}^\lambda(R)L_{\overline{1\rightarrow n}}\,)\\
&=&\hspace*{-2.5mm}\delta_{ik} \Tr_{R(1\dots n)}(e_{ii}^\lambda(R) L_{\overline{1\rightarrow n}}\,) = \delta_{ik}\,s_\lambda(L).
\end{eqnarray*}

So, we come to the final result:
$$
\Tr_{R(1\dots n)}\left(\rho_R(z)\,L_{\ov{1\to n}}\,\right)  = \sum_{\lambda\,\vdash n}s_{\lambda}(L)\sum_{i=1}^{d_\lambda}Z_{ii}^\lambda
= \sum_{\lambda\,\vdash n}s_{\lambda}(L)\,\chi^\lambda(z).
$$

To complete the proof we take $z=z_\nu$ of cyclic type $\nu$ and use the Poposition \ref{prop:11}. Observe that the value $\chi^\la_\nu=\chi(z_\nu)$
is the same for all elements $z_\nu$ of the same cyclic type $\nu$ since any two such elements are conjugated in the braid group. \hfill\rule{6.5pt}{6.5pt}

\medskip

Thus, we have established a $q$-version of the Frobenius formula by presenting the symmetric polynomials as elements of the RE algebra $\LL(R)$. However, it does not explain the appearance of the Hall-Littlewood polynomials (and their super-analogues) in this formula. To understand this fact we define the {\em eigenvalues} of the matrix $L$, which enables us to write down the symmetric polynomials associated with the RE algebras in another form. The key point in this definition is the quantum Cayley-Hamilton
identity for the generating matrix $L$.

For the sake of simplicity we first consider tthe RE algebras $\LL(R)$ associated with the even Hecke symmetries $R$ of bi-rank $(m|0)$. In this case the generating matrix $L$
satisfies the following Cayley-Hamilton identity:
\be
L^m-q\, e_1(L)\, L^{m-1}+q^2\, e_2(L)\, L^{m-2}+...+(-q)^{m-1}\,e_{m-1}(L)\, L+ (-q)^{m}\,e_{m}(L)\, I=0,
\label{CH}
\ee
where $e_k(L)$ are the quantum elementary symmetric polynomials. They are particular cases of the Schur polynomials $s_\la(L)$, corresponding to the partitions $\la=(1^k)$.
According to the general definition (\ref{schu}) these functions are defined by
$$
e_k(L)=\Tr_{R(1\dots k)}(\rho_R(a_k)L_{\ov{1\to k}}\,),
$$
where  $a_k \in\H_k(q)$ is the primitive orthogonal idempotent corresponding to one-column diagram $\la=(1^k)$. Note that in this case $d_\la=1$. The idempotents $a_k$ can
be defined by the following recursion:
$$
a_1 = e,\quad a_k=\frac{1}{k_q} \, a_{k-1}(q^{k-1} \, e-(k-1)_q\, \tau_{k-1})\, a_{k-1}, \quad k\ge 2.
$$

Observe that by a similar recursion it is possible to define the primitive idempotents $h_k$ parameterized by the partitions $\la=(k)$. The corresponding Schur polynomials are called
quantum full symmetric functions.

Now, we introduce  analogues of the eigenvalues of numerical matrices. Let $\mu_i$, $1\le i\le m$, be indeterminates, meeting the following system of algebraic equations:
\be
\sum_i^m \mu_i=q\, e_1(L),\qquad \sum_{i<j}^m \mu_i\, \mu_j=q^2 e_2(L),\quad \dots \quad  \prod_{i=1}^m \mu_i = q^{m}e_{m}(L).
\label{param}
\ee
We call $\mu_i$ {\em quantum eigenvalues} of the matrix $L$. We assume them to be central in the algebra $\LL(R)[\mu_1,\dots,\mu_m]$. It is interesting to parameterize all symmetric
polynomials of the algebra $\LL(R)$ via these eigenvalues. The elementary symmetric polynomials $e_i(L)$ are already parameterized according to formula (\ref{param}). In fact, these
polynomials are the usual elementary symmetric polynomials with respect to the quantities  $q^{-1}\mu_i$.  As for other Schur polynomials, they can be expressed via the elementary ones
with the use of the classical Jacobi-Trudi formulae (see \cite{GPS2} for detail).

The corresponding parametrization of the power sums $p_k(L)$ is as follows:
\be
p_k(L)=\sum_{i=1}^m\, \mu_i^k d_i,\qquad d_i=q^{-1}\, \prod_{j\not= i}^m \, \frac{\mu_i-q^{-2}\, \mu_j}{ \mu_i-\mu_j}.
\label{powss}
\ee
\begin{remark}\rm
From the parameterization (\ref{powss}) it is not evident that $p_k(L)$ are indeed polynomials in eigenvalues $\mu_i$. The simple way to verify this is to use  the set of quantum
Newton identities connecting $p_i(L)$ and $e_j(L)$:
$$
k_q e_k- q^{k-1}e_{k-1}\,p_1+q^{k-2}e_{k-2}\,p_2+\dots+(-1)^kp_k =0,\qquad \forall \,k\ge 1.
$$
The set of these identities allows one to express $p_k(L)$ as polynomials in $e_i(L)$ and vice versa. Now, by using the  parameterization (\ref{param}), we get the corresponding
polynomial parameterization of the quantum power sums.
\end{remark}

The polynomials $p_k(L)$, $k\ge 1$ are equal (up to a numerical factor $(q-q^{-1})$) to  the Hall-Littlewood polynomials, corresponding to the partitions $\la=(k)$ if we put
 $t=q^{-2}$, where $t$ is  the parameter entering the definition of the Hall-Littlewood polynomials in the notation of \cite{M}. Thus, the Hall-Littlewood polynomials are
 the quantum power sums (up to the indicated modifications).

Now, let us assume a skew-invertible Hecke symmetry $R$ to be of a general bi-rank $(m|n)$ with $m\,n\not=0$. In this case the generating matrix $L$ of the corresponding algebra
$\LL(R)$ satisfies a more complicated identity that that (\ref{CH}). Without going into details, which can be found in  \cite{GPS1, GPS2}, we only point out the main differences between
the general and even cases.

First, the degree of the corresponding Cayley-Hamilton polynomial is  $m+n$.  Second, this polynomial is not monic. Being multiplied by the quantum Schur function $s_\lambda(L)$
corresponding to the partition $\la=(n^m)$, it can be factorized in the product of two polynomials with central coefficients. One of them is of degree $m$, the second one is of degree
$n$. In a similar way as above it is possible to introduce the roots of these polynomial factors. Let $\mu_1,\dots, \mu_m$ and $\nu_1,\dots, \nu_n$ be the roots of the first and the
second polynomials respectively. We call all these roots the {\em eigenvalues} of the matrix $L$.  They are assumed to be central in the algebra $\LL(R)[\mu_1,\dots,\mu_m,\nu_1,\dots,
\nu_n] $. The eigenvalues $\mu_i$ are called  {\em even}, these $\nu_j$ --- {\em odd}.

Thus, we have two families  of eigenvalues --- even and odd. In \cite{GPS2} some quantum Schur polynomials were explicitly parameterized via the eigenvalues $\mu=\{\mu_i\}$
and $\nu=\{\nu_j \}$. In particular, we have:
$$
s_{(k)}(\mu, \nu)= \sum_{r=0}^k e_r(-q\nu)\,h_{k-r}(q^{-1}\mu),
$$
$$
s_{(1^k)}(\mu, \nu) =  \sum_{r=0}^k\,e_r(q^{-1}\mu)\, h_{k-r}(-q \nu),
$$
where $e_k$ and $h_k$ are the usual elementary and full symmetric polynomials respectively. As for the other quantum Schur polynomials, their parametrization via $\mu$ and $\nu$
can be also obtained with the use of Jacobi-Trudi formulae.

A spectral parametrization for the power sums, obtained in  \cite{GPS4}, is as follows:
$$
p_k(L)=\sum_{i=1}^m\, \mu_i^k\,d_i+\sum_{j=1}^n\, \nu_j^k \, {\tilde d_j},
$$
where
$$
d_i=q^{-1}\, \prod_{p\not= i}^m \, \frac{\mu_i-q^{-2}\, \mu_p}{ \mu_i-\mu_p}\, \prod_{j=1}^n \frac{\mu_i-q^{2}\, \nu_j}{ \mu_i-\nu_j},\qquad
{\tilde d_j}=-q\, \prod_{i= 1}^m \, \frac{\nu_j-q^{-2}\, \mu_i}{ \nu_j-\mu_i}\, \prod_{p\not=j}^n \frac{\nu_j-q^{2}\, \nu_p}{ \nu_j-\nu_p}.$$

Note, that all these polynomials are super-symmetric in the sets of variables $x=q^{-1} \mu$ and $y=q\, \nu$ in the following sense (see for instance \cite{S}). A polynomial $P(x,y)$ is
super-symmetric, if it is symmetric with respect to permutations of the  variables $x$ and $y$  separately and if we set $x_i=y_j =s$ for any fixed $i$ and $j$ the
polynomial $P$ does not depend on $s$.

The fact that the quantum Schur polynomials are super-symmetric is clear. A similar pro\-per\-ty for the power sums can be deduced from the Newton identities, which are also valid in
the general case.

Observe that for $q=1$ the power sums acquire the super-classical form:
$$
p_k(L)=\sum_{i=1}^m \mu_i^k-\sum_{j=1}^n \nu_j^k.
$$
It should be emphasized that this formula also takes place in the algebra $\LL(R)$ corresponding to any skew-invertible involutive symmetry $R$.

\begin{remark} \rm
As we noticed above, the quantum symmetric polynomials defined in the quantum matrix algebras related to pair of braidings $(R,F)$ (see \cite{IOP}) are not central if $F\not=R$.
In particular, this is true for the coefficients of the Cayley-Hamilton  polynomial, which belong to a commutative characteristic subalgebra, but they are not central. Consequently, the
centrality of the eigenvalues of the generating matrix $L$ also fails. Besides, the powers of the matrix $L$ entering the Cayley-Hamilton identity are not the usual matrix products.
The so-called RTT algebras are associated with the pair  $(R, P)$, where $P$ is the usual flip. One of the main discrepancy between RE algebras and the RTT ones in the case
when $R$ comes from the quantum group $U_q(sl(N))$ consists in the following. The RTT algebra cannot be equipped with the adjoint action of $U_q(sl(N))$, whereas for the RE
algebras such an action is well defined. Moreover, the map $L\mapsto L^k$, where $L$ is a generating matrix of a given RE algebra, is $U_q(sl(N))$-covariant. This is the reason,
why the Cayley-Hamilton identity for the generating matrix of the RE algebra is similar to the classical one, whereas for the RTT algebras it is not so.
\end{remark}

We conclude the paper by considering two examples, associated with the Hecke algebra $\H_3(q)$ and the following Hecke symmetries
$$
R^{(2)} =\left(\begin{array}{cccc}
q&0&0&0\\
0&q-q^{-1}&1&0\\
0&1&0&0\\
0&0&0&q
\end{array}\right),\qquad R^{(1|1)} =
\left(\begin{array}{cccc}
q&0&0&0\\
0&q-q^{-1}&1&0\\
0&1&0&0\\
0&0&0&-q^{-1}
\end{array}\!\!\right).
$$

We have 3 partitions of the number 3: $\la(1)=(3)$, $\la(2)=(2,1)$, $\la(3)=(1^3)$.  In this case the $q$-Frobenius formula (\ref{q-Fr}) takes the form:
\be
\left(\!\begin{array}{l}
p_{(3)}\\
p_{(2,1)}\\
p_{(1^3)}
\end{array}\!\!\!\right)=\left(\begin{array}{lcr}
q^2&-1&q^{-2}\\
q&q-q^{-1}&-q^{-1}\\
1&2&1\phantom{qq}
\end{array}\right)\left(\!
\begin{array}{l}
s_{(3)}\\
s_{(2,1)}\\
s_{(1^3)}
\end{array}\!\!\!\right).
\label{q-Fr3}
\ee
The $3\times 3$ matrix of $\H_3(q)$ irreducible characters in this formula slightly differs from that presented in \cite{R} since we use another normalization for the Hecke algebra
generators and the parameter $q$.

Now, we write explicitly all needed quantum symmetric polynomials in terms of eigenvalues of the generating matrix $L$ in the RE algebra $\LL(R)$ corresponding to the first Hecke symmetry $R^{(2)}$. The matrix $R^{(2)}$ comes from the quantum group $U_q(sl(2))$ and in this case $L$ has two even eigenvalues $\mu_1$ and $\mu_2$. The quantum Schur
polynomials and power sums have the following explicit forms:
$$
\begin{array}{lcl}
p_1=q^{-1}(\mu_1+\mu_2) & & s_{(3)}=q^{-3}(\mu_1^3+\mu_2^3+\mu_1\, \mu_2\,(\mu_1+ \mu_2)) \\
\rule{0pt}{5mm}
 p_{2}=q^{-1} \mu_1^2+ \mu_2^2+q^{-2}(q-q^{-1})\mu_1\, \mu_2 & & s_{(2,1)}=q^{-3} \mu_1\, \mu_2\,(\mu_1+ \mu_2) \\
\rule{0pt}{5mm}
p_{3}=q^{-1} (\mu_1^3+ \mu_2^3)+q^{-2}(q-q^{-1})\mu_1\, \mu_2\, (\mu_1+\mu_2)  & & s_{(1^3)}=0
\end{array}
$$
The second Hecke symmetry $R^{(1|1)}$ comes from the quantum super-group $U_q(sl(1|1))$. The generating matrix $L$ of the corresponding RE algebra $\LL(R)$ has one even
and one odd eigenvalue $\mu$ and $\nu$ respectively. The quantum super-symmetric polynomials are as follows:
$$
\begin{array}{lcl}
p_1=q^{-1}\mu-q\,\nu & \qquad& s_{(3)}=q^{-2}\mu^2\,p_1\\
\rule{0pt}{5mm}
p_{2}=(\mu+ \nu)\,p_1 & & s_{(2,1)}=-\mu\,\nu \,p_1\\
\rule{0pt}{5mm}
p_{3}=(\mu^2+ \mu\nu+\nu^2)\,p_1 & & s_{(1^3)}=q^{2} \nu^2\, p_1.
\end{array}
$$
Taking into account that $p_{(3)} = p_3$, $p_{(2,1)} = p_2\,p_1$ and $p_{(1^3)} = p_1^3$ it is not difficult to verify that the $q$-Frobenius formula (\ref{q-Fr3}) is valid in the both cases.

\bigskip

\leftline{\bf Data Availability}

\medskip

\noindent
Data sharing not applicable to this article as no datasets were generated or analysed during the current study.

\end{document}